\pdfoutput=1 % uncomment to ensure pdflatex processing (mandatory e.g. to submit to arXiv)

\documentclass[justified,notoc,twoside]{tufte-handout}
\usepackage[T1]{fontenc}
\usepackage[english]{babel}
\usepackage[tracking]{microtype}
\usepackage{listings}
\usepackage[inline]{asymptote}

\usepackage[strict]{changepage}
\usepackage{amsthm}
\usepackage{thmtools}
\usepackage{thm-restate}
\usepackage{hyperref}
\usepackage{accsupp} % suppress sidenote marks from copy-paste
\usepackage{cleveref}

\usepackage{beton}
\usepackage{sectsty}
\allsectionsfont{\fontfamily{LinuxBiolinumT-OsF}\selectfont}

\usepackage{amsmath} % extended mathematics
\usepackage{amssymb}
\usepackage{units} % non-stacked fractions and better unit spacing
\usepackage[osf,sc]{mathpazo} % With old-style figures and real smallcaps.

\usepackage[most]{tcolorbox} % For rounded corner boxes
\usepackage{tocloft}   % Optional: For better ToC customization
\cftsetindents{section}{0em}{2em}
\cftsetindents{subsection}{0em}{2em}
\let\subsubsection\subsection

\usepackage[caption=false]{subfig}
\usepackage{graphicx} % allow embedded images
\setkeys{Gin}{width=\linewidth,totalheight=\textheight,keepaspectratio}
\graphicspath{{figures/}} % set of paths to search for images
\usepackage{booktabs} % book-quality tables
\usepackage{multicol} % multiple column layout facilities
\usepackage{fancyvrb} % extended verbatim environments
\fvset{fontsize=\normalsize} % default font size for fancy-verbatim environments

\newcommand{\bothcite}[2][0cm]{\citep{#2}\cite[#1]{#2}}

\AtBeginDocument{ \fontsize{12}{14}\selectfont } % increase font size

\title[GP II.1] % short margin title should not be more than half width of page
{Graph Puzzles II.1:\\
Counterexamples to\\ Jain's second unit vector flows conjecture} % [header short title] {Full title}

\date{}

\begin{document}

\selectlanguage{english}

\maketitle % this prints the handout title, author, and date

%%%%%%%%%%% Start Author Section %%%%%%%%%%%%%%%%%%
Nikolay Ulyanov\textsuperscript{*}\marginnote{\textsuperscript{*}{\scriptsize ulyanick@gmail.com}}
%%%%%%%%%%% End Author Section %%%%%%%%%%%%%%%%%%

\begin{abstract}

\noindent

\newthought{Abstract: }A $3$-dimensional nowhere-zero flow on a graph $G$ is a flow where each edge is assigned a $3$-dimensional vector with unit norm (which corresponds to the points of a $2$-dimensional unit sphere $S^2$). K.~Jain posed two conjectures related to this idea. First one suggests that such a flow exists for all bridgeless graphs. The~second conjecture states that we can assign values $\{-4,-3,-2,-1,1,2,3,4\}$ to the points of~$S^2$, such that antipodal points get opposite values, and values of any three equidistant points on great circles sum to zero. If both conjectures would be true, together they would imply Tutte's 5-flow conjecture. We show 2 counterexamples to the second conjecture, by constructing sets of points each of which additionally requires values $\{-5, 5\}$.
\end{abstract}

\noindent\rule{5in}{0.4pt}

\begin{figure}
  \checkoddpage \ifoddpage \forcerectofloat \else \forceversofloat \fi
  \hspace*{-1.5cm}
  \includegraphics[width=1.2\linewidth]{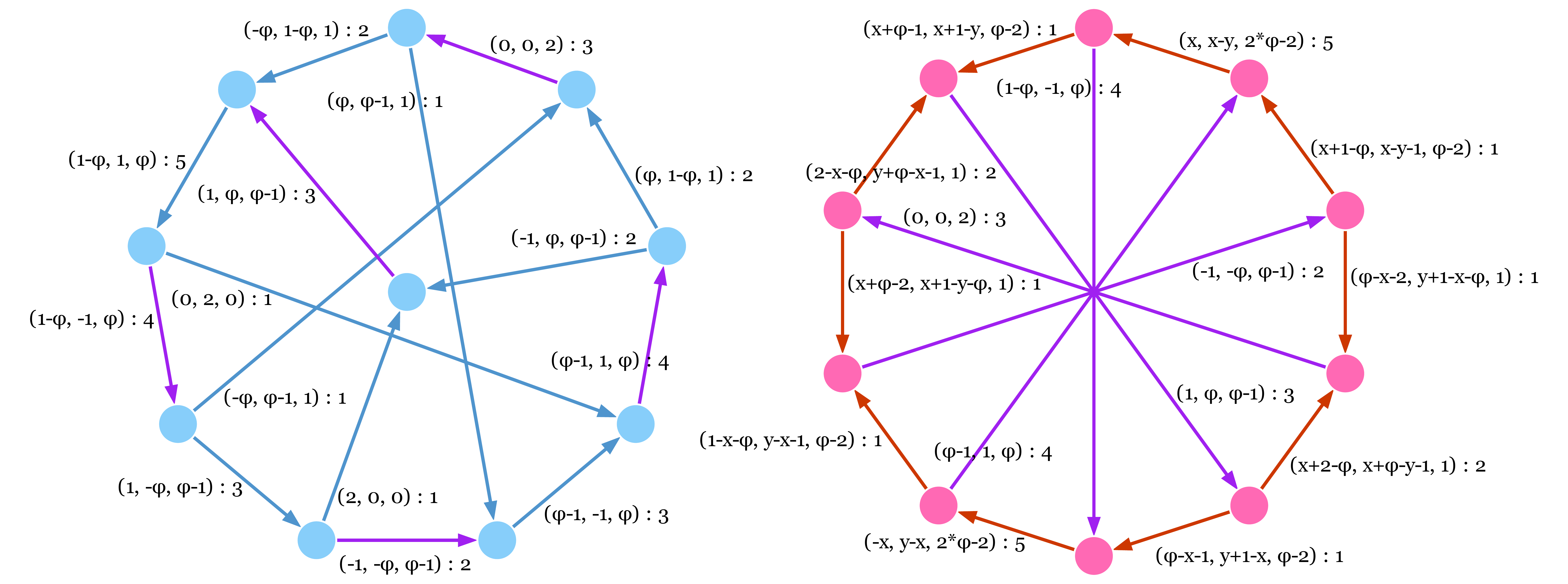}
  \caption{First counterexample to Jain's second conjecture, which might also be reinterpreted as a combination of Petersen graph and M\"obius ladder; scaled to radius~2}
  \label{fig:counterexample1_petersen_moebius}
  \setfloatalignment{b}
\end{figure}

\pagebreak

\begin{tcolorbox}[
    breakable,
    toggle enlargement=evenpage,
    rounded corners, % Rounded corners
    boxrule=1pt, % Thickness of the border
    width=1.2\linewidth, % Box width
]
    \csname @starttoc\endcsname{toc}
\end{tcolorbox}

\section{Introduction}\label{sec:Introduction}

The Open Problem Garden page \bothcite{OPGUVF} for unit vector flows states 2 conjectures, both by K.~Jain:

\begin{restatable}[unit vector flows]{conj}{}
\label{conj1}
Every bridgeless graph has a flow $f: E(G) \rightarrow S^2 = \left\{x \in \mathbb{R}^3 : |x| = 1\right\}.$
\end{restatable}

\begin{restatable}[$S^2$ nz5-flow]{conj}{}
\label{conj2}
There exists a map $q : S^2 \rightarrow \{-4,-3,-2,-1,1,2,3,4\}$ so that antipodal points sum to 0, and any three points which are equidistant on a great circle also sum to 0.
\end{restatable}

Would both conjectures be true, together they would imply the famous conjecture by W.~T.~Tutte \bothcite{Tutte}:

\begin{restatable}[nz5]{conj}{}
Every bridgeless cubic graph has a nowhere-zero 5-flow.
\end{restatable}

A \textit{nowhere-zero flow} ($nz$-flow) is a non-zero valuation on the edges $f(e), e \in E$ such that the sum of incoming values equals the sum of outgoing values. If all values are non-zero integers such that $f(e) < k$, then we say that it's a \textit{nowhere-zero $k$-flow} ($nzk$-flow).

The second conjecture has recently seen increase in attention, e.~g. it appeared in \bothcite{Li2026} and \bothcite{Houdrouge2026}. In this paper we will disprove it, by finding two $S^2$ point subsets, each of which requires values $\{\pm1,\pm2,\pm3,\pm4,\pm5\}$. But as an icebreaker, first we consider a non-counterexample of a nz5-flow on the sphere.

\section{Icosidodecahedron, Petersen graph and $nz5$-flow}

One of the first ideas that comes to mind when searching for $S^2$-flows, is to check various famous point configurations on the sphere, e. g., vertices of regular polyhedra (Platonic solids), Catalan solids, Archimedean solids, etc. One quickly finds, the following example: an Archimedean solid known as icosidodecahedron, which has 30 points (or 15 pairs of opposite points). Let $\varphi$ denote the golden ratio $\frac{1 + \sqrt{5}}{2}$. Then we can consider all even permutations of the following points as points of icosidodecahedron $\mathcal{ICOS}$:

$$
\left(0, 0, \pm 1\right), \left(\pm \frac{\varphi - 1}{2}, \pm \frac{1}{2}, \pm \frac{\varphi}{2}\right).
$$

As an example, here are a couple of triples of points which sum up to zero:

$$
(0, 0, 1), \left(\frac{\varphi}{2}, \frac{\varphi - 1}{2}, -\frac{1}{2}, \right), \left(-\frac{\varphi}{2}, -\frac{\varphi - 1}{2}, -\frac{1}{2}, \right),
$$

and

$$
\left(\frac{\varphi}{2}, \frac{\varphi - 1}{2}, \frac{1}{2}, \right), \left(-\frac{\varphi - 1}{2}, \frac{1}{2}, -\frac{\varphi}{2}\right), \left(-\frac{1}{2}, -\frac{\varphi}{2}, \frac{\varphi - 1}{2}\right).
$$

\begin{figure}
  \checkoddpage \ifoddpage \forcerectofloat \else \forceversofloat \fi
  \includegraphics[width=0.93\linewidth]{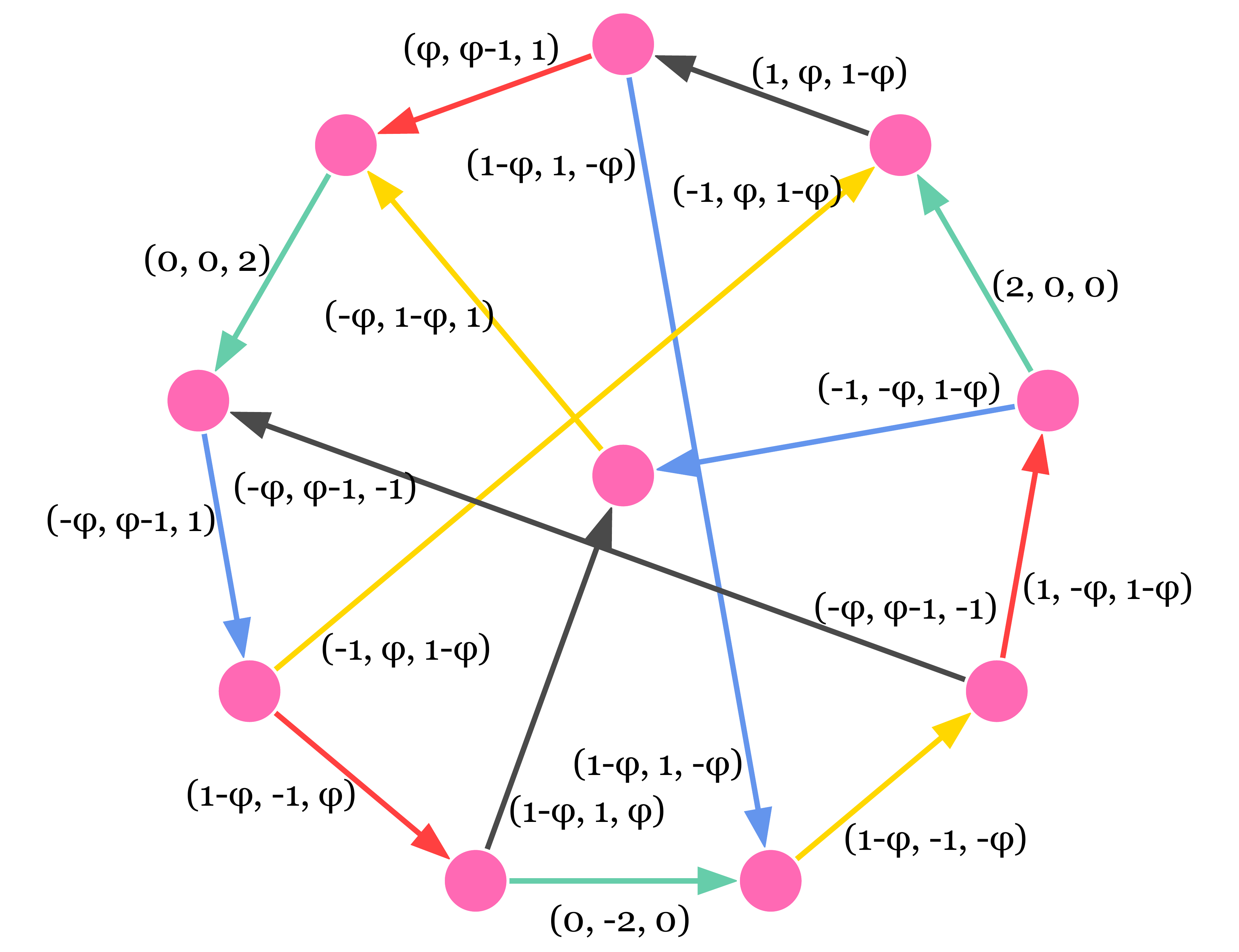}
  \caption{Icosidodecahedron mapping for Petersen graph, scaled to radius~2}
  \label{fig:petersen_icosi_flow}
  \setfloatalignment{b}
\end{figure}

There are 30 points, and 20 triples in total, each point appears exactly in 2 triples. If we consider opposite points to be equivalent, then we have 15 pairs of points, 10 pairs of triples, which we can reinterpret as a cubic graph. Specifically we get the Petersen graph (see Fig.~\ref{fig:petersen_icosi_flow}), which is also the smallest snark, and so it doesn't have a $nz4$-flow. It has various $nz5$-flows, so we also can construct an example of a map $\mathcal{ICOS} \rightarrow \{-4,-3,-2,-1,1,2,3,4\}$.

The same construction has also appeared in \bothcite{Mattiolo2025} and \citep{Houdrouge2026}, although without explicit coordinates and identification with icosidodecahedron. Also interesting to note that the 10~pairs of triples mentioned above, each of them lives on some great circle, but they are not immediately visible on icosidodecahedron, the more direct correspondence would be with its relative non-convex polyhedron named dodecadodecahedron.

\section{Counterexamples to $S^2$ nz5-flow conjecture}\label{sec:counterexamples}

We now present two counterexamples to Conjecture~\ref{conj2}. Both are finite subsets of $S^2$ for which no labeling ${q : S^2 \rightarrow \{-4,-3,-2,-1,1,2,3,4\}}$ satisfying the conjecture's constraints exists. In each case, the labeling requires values $\pm 5$, i.e., a nowhere-zero $6$-flow suffices but a nowhere-zero $5$-flow does not.

\subsection{Counterexample 1: 50-point expansion of the icosidodecahedron}\label{sec:ce1}

\subsubsection*{Construction}

The first counterexample is obtained by expanding the icosidodecahedron 30 point set. For each point $p$ of icosidodecahedron we consider a corresponding small circle, equidistant from the point $p$, which lies on a plane $\langle x, p \rangle = -0.5$ (assuming the sphere has radius 1 and centered at the origin). Next we consider all pairs of points of icosidodecahedron that lie at a spherical distance of $2\pi/5$, which is double the angular distance between adjacent vertices on the 6 great circles of icosidodecahedron (they are different from the 10 great circles mentioned above). For each such pair of points we intersect the two corresponding small circles. This produces new points on~$S^2$ that, together with the original 30, form a set of 50 points with 40 great-circle triples.

\subsubsection*{Coordinates}

Let $\varphi = \frac{1+\sqrt{5}}{2}$ denote the golden ratio again, and let $x = \frac{2}{\sqrt[\leftroot{-2}\uproot{2}4]{5}
}$, $y = x\varphi$. Scaling the sphere by $R = 2$, the 50 points have coordinates that are integer linear combinations of $1, \varphi, x, y$. The full list of 50 coordinates and 40 triples is given in the accompanying code repository and in Fig.~\ref{fig:counterexample1_petersen_moebius}.

If we quotient by antipodes and regard each triple as a vertex, then besides the Petersen graph we can recover another cubic graph\,---\,a $10$-vertex M\"obius ladder. The 25 point-pairs split into 10 Petersen-only old edge-orbits, 10 new edge-orbits coming from the added points, and 5 old edge-orbits shared by both quotients, the latter of which correspond to perfect matchings in both graphs. An explicit $nz6$ labeling of these 25 point-pairs can be seen in Fig.~\ref{fig:counterexample1_petersen_moebius}.

\subsubsection*{Verification}

We encode the labeling problem as a Boolean satisfiability (SAT) instance: for each of the 25 antipodal representatives, we introduce 8 Boolean variables (one per allowed value in $\{-4,\ldots,-1,1,\ldots,4\}$), with exactly-one constraints, antipodal constraints ($q(-p) = -q(p)$), and triple-sum constraints ($q(a)+q(b)+q(c)=0$). The resulting CNF formula (200 variables, 19765 clauses) is \textbf{unsatisfiable}, as confirmed by the PicoSAT/pycosat solver. This proves that no labeling with values in $\{-4,\ldots,-1,1,\ldots,4\}$ exists for this 50-point set.\sidenote[][-2.5cm]{We couldn't escape the recent trend and also verified the construction in Lean~4, even though it doesn't make much sense in our case and doesn't add anything new. Nevertheless, using the \texttt{bv\_decide} tactic (which is basically a SAT solver), we also provide a formally verified proof in the accompanying GitHub repository.}

\subsection{Counterexample 2: 36-point construction based on square root arithmetic}\label{sec:ce2}

\subsubsection*{Construction}

The second counterexample uses a slightly different construction, which produces a smaller 36-point subset of $S^2$ that also requires a $nz6$-flow.

We fix parameters $v_1 = 1$, $v_2 = 3$, and $w = 2$. We generate a set of candidate coordinate values as follows. For integers $w_1 \in \{0, 1, \ldots, w\}$ and $w_2 \in \{-w, \ldots, w\}$, we compute two types of coordinates:
$$
c_{\text{rat}} = \frac{|w_1 \sqrt{v_1} + w_2\sqrt{v_2}|}{2}, \qquad c_{\text{sqrt}} = \sqrt{\frac{|w_1\sqrt{v_1} + w_2\sqrt{v_2}|}{2}}.
$$
We only keep values $c \le 1$. We then search for all triples $(c_1, c_2, c_3)$ from this combined set satisfying $c_1^2 + c_2^2 + c_3^2 = 1$ (i.e., points on $S^2$). Each such triple generates up to 48 points on $S^2$ via all permutations and sign changes of coordinates.

We identify all great-circle triples, and after iteratively removing triples where at least 2 of the 3 vertices appear in only one triple, we obtain a connected component of 126 points and 108 triples. After additional pruning we finally converge to a subset of 36 points and 13 triples.

\subsubsection*{Coordinates}

The 7 distinct absolute coordinate values are:
$$
0,\; \frac{2-\sqrt{3}}{2},\; \frac{\sqrt{3}-1}{2},\; \frac{1}{2},\; \sqrt{\sqrt{3}-1},\; \frac{\sqrt{3}}{2},\; 1.
$$

\subsubsection*{Verification}

The impossibility of a $nz5$-flow labeling is verified by SAT solving: the CNF formula (144 variables, 6710 clauses) is \textbf{unsatisfiable}, confirmed by PicoSAT/pycosat. A labeling with values in $\{-5,\ldots,-1,1,\ldots,5\}$ (a $nz6$-flow) does exist, found by SAT solving with 180 variables and 13048 clauses.

\section{Final remarks}

Both counterexamples disprove Conjecture~\ref{conj2} by exhibiting finite point sets on $S^2$ where the labeling constraints force the use of values $\pm 5$. Several natural questions arise:

\begin{itemize}
\item Can we find a subset which additionally requires the value $\pm 6$ or more? Thus far our best attempts only produced examples where $\pm 5$ is needed at maximum.
\item What is the smallest subset that still works as a counterexample to Conjecture~\ref{conj2}?
\item Is it still possible to salvage the Conjecture~\ref{conj2} for the purposes of combining it with Conjecture~\ref{conj1} to prove Tutte's 5-flow conjecture? One could imagine here, that we would need to find a (possibly non-finite) subset of points on $S^2$, maybe with some nice properties (e. g. its points should have algebraic coordinates of some specific form), which still admits a nowhere-zero 5-flow mapping, and is suitable for reuse in Conjecture~\ref{conj1} (so we can map all cubic bridgeless graphs to this subset). We have done some research in this direction, which will be shared in the next paper in the series \bothcite{GP2-2}.
\item In all our experiments, we never observed a difference between nowhere-zero $k$-flows and nowhere-zero mod-$k$ flows, in the sense that any subset we checked has the same minimal value of $k$ for both kinds of flows, and it would be interesting to understand whether this is always the case or not.
\end{itemize}

 We also note that numerical precision is important in these computations: we use tolerance $\varepsilon = 10^{-7}$ for point identification and great-circle triple detection, and all results have been cross-validated with exact arithmetic where possible.

The code for both counterexamples, including the SAT encoding, backtracking search, and Lean~4 verification, is available at the accompanying GitHub repository.\sidenote[][-0.5cm]{\url{https://github.com/gexahedron/unit-vector-flows}}.

\bibliography{gp2-1}

\bibliographystyle{alpha}

\end{document}